\documentclass[10pt]{article}
\usepackage[a4paper,hmargin=29mm,vmargin=38mm]{geometry}

\usepackage[utf8]{inputenc}
\usepackage[T2A,T1]{fontenc}
\usepackage{lmodern}
\usepackage{fix-cm}
\usepackage{textcomp}

\usepackage[russian,french]{babel} 
\usepackage[french]{varioref}
\usepackage[autolanguage]{numprint}

\usepackage[hyphens]{url}%
\RequirePackage{etoolbox}
\appto\UrlSpecials{%
  \do\R{\penalty0 \mathchar`\R }%
  \do\U{\penalty0 \mathchar`\U }%
  \do\D{\penalty0 \mathchar`\D }%
  \do\F{\penalty0 \mathchar`\F }%
  \do\A{\penalty0 \mathchar`\A }%
  \do\W{\penalty0 \mathchar`\W }%
  \do\Z{\penalty0 \mathchar`\Z }%
  \do\G{\penalty0 \mathchar`\G }%
  \do\P{\penalty0 \mathchar`\P }%
  \do\N{\penalty0 \mathchar`\N }%
}
\providecommand\url[1]{\texttt{\detokenize{#1}}}
\usepackage{csquotes}
\usepackage[numbers,comma,square]{natbib}
\setcitestyle{square,comma,numbers}
\let\hyphenate\-

\usepackage{amsmath}
\usepackage{amsfonts}
\usepackage{mathtools} \mathtoolsset{mathic} % kerning of math surrounded by italic

\usepackage{amsthm} % ams theorems that are not lists
\theoremstyle{plain}
\newtheorem*{trente}{Proposition~30}
\newtheorem*{vingt}{Proposition~20}
\newtheorem*{treize}{Proposition~13}
\newtheorem*{dixneuf}{Proposition~19}
\newtheorem*{six}{Proposition~6}
\newtheorem*{porisme}{Porisme}
\newtheorem*{idb}{Invariants de boucle}
\newtheorem*{dvingt}{Définition~20}
\newtheorem*{fix}{La réparation}
\theoremstyle{definition}
\newtheorem{exemple}{Exemple}
\newtheorem*{remerciements}{Remerciements}

\makeatletter
\renewcommand\section{\@startsection{section}{1}{0pt}{\medskipamount}{-.5em}{\normalfont\normalsize\bfseries}}
\makeatother
\usepackage{titlesec}
\titlelabel{\thetitle.\ }

\usepackage{indentfirst}
\ifx\texorpdfstring\undefined\newcommand\texorpdfstring[2]{#1}\fi
\ifx\pdfbookmark[3]{}\fi

\usepackage{color}
\usepackage{array}
\usepackage{hhline}
\DeclareMathOperator\pgcd{pgcd}

\begin{document}
\title{Euclide avait-il besoin de l'algorithme d'Euclide pour démontrer l'unicité de la factorisation?}
\author{David Pengelley \and Fred Richman}
\date{}
\maketitle

Cet article est la traduction française, par Henri Lombardi et Stefan Neuwirth, de l'article ``Did Euclid need the Euclidean algorithm to prove unique factorization?'' écrit par David Pengelley et Fred Richman. Il s'agit d'un article paru à l'origine dans l'American Mathematical Monthly \no113 en~2006, pages~196-205. 

Nous remercions les auteurs et la Mathematical Association of America (analogue de l'APMEP en France) de nous autoriser la publication de cette traduction française.

\section{Introduction.}
 Le théorème fondamental de l'arithmétique affirme que tout entier naturel est, de manière unique, un produit de nombres premiers. Le cœur de cette unicité se trouve dans le Livre~VII des \emph{Éléments} d'Euclide~\cite{vitrac94}:
\begin{trente}[lemme d'Euclide]
  Si un nombre premier divise un produit, il divise l'un des facteurs.
\end{trente}

Euclide commence le Livre~VII en introduisant l'algorithme d'Euclide. À partir de sa preuve que l'algorithme fonctionne, il déduit un résultat algébrique: 

\begin{porisme}[propriété algébrique du pgcd]
Si un nombre divise deux nombres, il divise leur plus grand commun diviseur.\footnote{N.D.T. Le mot «porisme» utilisé dans les \emph{Éléments} est aujourd'hui inusité. Une traduction moderne serait «corollaire immédiat». Ce porisme est donc un \emph{porisme de l'algorithme d'Euclide}.}
\end{porisme}
\noindent Le lemme d'Euclide peut être déduit à partir de cette propriété, mais il n'est pas vraiment évident qu'Euclide le fasse. Nous serions très surpris qu'Euclide n'utilise pas cette propriété, car il la signale très tôt et parce que nous nous attendons à ce qu'il fasse usage de l'algorithme d'Euclide d'une manière significative. Dans cet article, nous explorons la question de savoir à quel point la propriété algébrique du pgcd intervient dans le lemme d'Euclide, si c'est effectivement le cas. 

Une idée centrale dans la rédaction d'Euclide est celle de la proportionnalité entre quatre nombres: $a$~est à~$b$ comme $c$~est à~$d$. Euclide donne deux définitions différentes de la proportionnalité, une dans le Livre~VII pour les nombres, la «proportionnalité pythagoricienne», et une dans le Livre~V pour les grandeurs générales, la «proportionnalité eudoxienne». Nous allons découvrir qu'il est essentiel de garder à l'esprit la distinction entre ces deux définitions, et que de nombreuses autorités, peut-être y compris Euclide lui-même, sont tombées dans le piège de croire que l'on peut voir facilement que la proportionnalité eudoxienne est la même que la proportionnalité pythagoricienne.

Pour terminer, nous suggérerons un moyen d'amender la démonstration d'Euclide après 2300~ans.

\section{L'algorithme d'Euclide.}

La théorie des nombres d'Euclide est contenue dans les livres VII~à~IX des \emph{Éléments}. Au début du livre~VII, il présente l'algorithme d'Euclide. L'entrée de l'algorithme est une paire de nombres (entiers strictement positifs) $a$~et~$b$ avec $a<b$, et l'algorithme consiste en une répétition indéfinie de trois étapes:
\[\text{Répéter }\left\{
  \begin{aligned}
1.&\text{ si $a$~divise~$b$ retourner~$a$;}\\
2.&\text{ tant que $a<b$ poser $b=b-a$;}\\
3.&\text{ poser $(a,b)=(b,a)$.}
  \end{aligned}\right.
\]
L'étape~2 est l'\emph{algorithme de division}: nous soustrayons~$a$ à~$b$ jusqu'à ce que $b$~devienne plus petit que~$a$. Nous pouvons aussi écrire $b=qa+r$, où $r<a$, et remplacer~$b$ par~$r$. À l'étape~3, on échange les rôles de $a$~et~$b$, car $b$~est maintenant le plus petit des deux.

L'algorithme d'Euclide est réputé retourner le plus grand commun diviseur de~$a$ et~$b$. Pour assurer ce résultat, il faut une démonstration, que fournit Euclide. Sa démonstration est essentiellement la première partie du théorème suivant, dont nous laissons la vérification au lecteur.
\begin{idb}
Si \(b=qa+r\), les affirmations suivantes sont vraies:
\begin{enumerate}
\item les diviseurs communs à \(a\)~et~\(b\) sont exactement les diviseurs communs à \(a\)~et~\(r\);
\item le sous-groupe des entiers engendré par \(a\)~et~\(b\) est égal au sous-groupe engendré par \(a\)~et~\(r\).
\end{enumerate}
\end{idb}
\noindent En appliquant ce résultat chaque fois que l'on effectue la boucle dans l'algorithme, nous voyons que l'ensemble des diviseurs communs à $a$~et~$b$ est un \emph{invariant de boucle~}: il est le même après avoir effectué les trois étapes que ce qu'il était auparavant. Ainsi, lorsque l'on sort de l'algorithme, ce qui arrive lorsque $a$~divise~$b$, nous avons la garantie que le plus grand commun diviseur est rendu parce que, lorsque~$a$ divise~$b$, $a$~est le plus grand commun diviseur de $a$~et~$b$. Euclide prouve d'abord que la sortie de l'algorithme est un diviseur commun de $a$~et~$b$, puis, pour démontrer que c'est \emph{le plus grand}, il montre que tout autre diviseur commun doit le diviser, et donc être plus petit. C'est cette propriété algébrique du pgcd qu'Euclide note dans le porisme. 

Cette propriété algébrique est le fait théorique révélé par l'analyse de l'algorithme d'Euclide. Un algorithme efficace qui calcule le plus grand commun diviseur ne sera pas un brise-glace théorique pertinent. Vous ne pourrez rien en déduire d'intéressant, car deux nombres quelconques admettent un plus grand commun diviseur pour la simple raison que l'ensemble des diviseurs communs est fini. En revanche, le fait que tout autre diviseur commun divise le plus grand commun diviseur est surprenant et riche de conséquences. C'est seulement à travers le porisme que l'algorithme d'Euclide peut jouer un rôle réel dans la théorie des nombres d'Euclide.

À la lumière du porisme, nous pouvons remplacer la notion de \emph{plus grand} commun diviseur par une notion purement algébrique: \emph{un pgcd algébrique} de deux nombres est un diviseur commun qui est divisible par tout autre diviseur commun. Il n'y a aucune raison de croire \emph{a priori} que deux nombres arbitraires possèdent un pgcd algébrique, mais c'est exactement ce que nous dit le porisme.

La deuxième partie du théorème nous dit que le sous-groupe engendré par $a$~et~$b$ est un invariant de boucle. À la fin, lorsque $a$~est égal au pgcd des $a$~et~$b$ d'origine, elle dit que $\pgcd(a,b)$~est dans le sous-groupe engendré par $a$~et~$b$, c'est-à-dire que l'on peut écrire
\[
\pgcd(a,b)=sa+tb
\]
pour certains entiers $s$~et~$t$. Ceci est connu comme l'\emph{équation de Bézout. }À partir de cette équation, il est facile de démontrer le porisme d'Euclide, que tout diviseur commun à $a$~et~$b$~\emph{doit} diviser~$\pgcd(a,b)$. 

De nos jours on prouve souvent le lemme d'Euclide en utilisant l'équation de Bézout. Supposons que $p$~est un nombre premier qui divise~$ab$, si $p$~ne divise pas~$a$, $p$~et~$a$ n'ont pas de diviseur commun non trivial, ainsi l'équation de Bézout dit qu'il existe des entiers $s$~et~$t$ tels
que 
\[
sp+ta=1\text.
\] 
En multipliant cette équation par~$b$ nous obtenons
\[
spb+tab=b\text,
\] 
ce qui montre que $p$~divise~$b$ puisqu'il divise~$ab$.

On peut également utiliser la propriété algébrique du pgcd pour démontrer que si $p$~ne divise pas~$a$, alors $\pgcd(pb,ab)=b$,\footnote{N.D.T. En effet $b$~doit diviser~$\pgcd(pb,ab)$ d'après le porisme; on écrit alors ce pgcd sous forme~$bu$, et l'on voit que $u$~est un diviseur commun à $p$~et~$a$, donc est nécessairement égal à~1 si $p$~ne divise pas~$a$.} ainsi $p$~divise~$b$. En explorant comment Euclide lui-même démontre le lemme d'Euclide, nous allons observer avec soin s'il fait appel à cette propriété.

\section{À la recherche d'une démonstration.}

Nous présentons maintenant une histoire légèrement mythologisée d'une quête d'une démonstration du lemme d'Euclide. Plus précisément, notre histoire commence avec la question: est-ce qu'Euclide a quelque chose d'intéressant à nous dire sur comment démontrer le lemme d'Euclide?

Comme indiqué précédemment, Euclide démarre le Livre~VII avec l'algorithme d'Euclide accompagné du porisme qui affirme que le plus grand commun diviseur de deux nombres est divisible par tout autre diviseur commun. Sa preuve du lemme d'Euclide réfère aux propositions 20~et~19.
\begin{vingt}
  Si \(u\)~et~\(v\) sont les plus petits entiers tels que \(u:v=c:d\), alors \(u\)~divise~\(c\) et \(v\)~divise~\(d\).
\end{vingt}
\noindent Nous allons déduire le lemme d'Euclide de la proposition~20 plus ou moins de la même manière que lui. Supposons qu'un nombre premier $p$~divise~$ab$, disons $ab=pc$. Considérons la fraction
\[
\frac ap=\frac cb\text.
\]
Si $u$~et~$v$ sont les plus petits nombres tels que $u/v=a/p$, la proposition~20 nous informe que $v$~divise à la fois $p$~et~$b$. Donc $v=1$ ou $v=p$, parce que $p$~est premier. Dans le premier cas $p$~divise~$a$, dans le second $p$~divise~$b$. Donc la proposition~20 est évidemment la proposition clé. Comment Euclide la démontre-t-il?

Nous suivons la paraphrase de Heath \cite{heath26} de la démonstration d'Euclide de la proposition~20. Euclide démontre que $u$~divise~$c$ (et donc $v$~divise~$d$), ou comme il l'écrit, que $u$~est une \emph{partie} de~$c$. Cela signifie que 
\[
u=\frac cn
\]
pour un certain entier positif~$n$. Il fait ceci en excluant l'alternative, à savoir que $u$~n'est pas une partie de~$c$, auquel cas nous pourrions écrire
\[
u=m\cdot \frac cn\text,
\]
où $n$~divise~$c$, et $m>1$. C'est-à-dire $c/n$~est une $n$ième partie de~$c$, et $u$~est égal à~$m$ de ces~$n$ièmes parties. Dans la proposition~4, Euclide a montré comment calculer de tels nombres $m$~et~$n$ à partir de $u$~et~$c$ en utilisant l'algorithme d'Euclide, et Heath réfère à la proposition~4 dans sa paraphrase.

Euclide affirme alors que $u:c=v:d$. Ceci est justifié par une autre proposition. 
\begin{treize}[alternando]
  Si \(a:b=c:d\), alors \(a:c=b:d\).
\end{treize}
\noindent C'est bon, car si $a/b=c/d$, alors $a/c=b/d$.

En poursuivant la démonstration de la proposition~20, Euclide note que
\[
v=m\cdot \frac dn\text.
\]
 Son énonciation pour ceci est «$v$~est les mêmes parties de~$d$ que~$u$ de~$c$», c'est-à-dire, $v$~est égal à $m$~$n$ièmes parties de~$d$ exactement comme $u$~est égal à $m$~$n$ièmes parties de~$c$. Parce que $m>1$, les nombres $c/n$~et~$d/n$ sont plus petits que $u$~et~$v$. Mais comme $c/n:c=d/n:d$, il s'ensuit que $c/n:d/n=c:d$, ce qui contredit le fait que $u$~et~$v$ seraient les plus petits nombres avec cette propriété. Ainsi $u$~doit diviser~$c$.

Il y a un sérieux problème avec ce raisonnement. Nous savons que $n$~divise~$c$, parce que c'est ainsi que nous avons choisi~$n$, mais pourquoi $n$~divise-t-il~$d$, c'est-à-dire, pourquoi $d/n$~est-il un nombre (entier)? Il est vrai que l'algorithme d'Euclide produit des nombres $m$~et~$n$ premiers entre eux (quoique Euclide ne mentionne pas ce fait), et nous savons que $nv=md$. Mais conclure de ces deux faits que $n$~divise~$d$ nécessite plus que le lemme d'Euclide lui-même, et c'est ce lemme que nous voulons finir par démontrer! Ceci a conduit Zeuthen \cite[pages~127-129]{zeuthen02} à dire que la démonstration par Euclide de son lemme est sans valeur, parce qu'Euclide devait supposer quelque chose d'essentiellement plus fort que le lemme lui-même pour démontrer la proposition~20. Ainsi, au lieu de trouver une démonstration alternative du théorème fondamental de l'arithmétique, nous avons trouvé une erreur chez Euclide!

Attendez une minute. C'est d'\emph{Euclide} que nous sommes en train de parler, de l'auteur du traité de mathématiques le plus fameux de tous les temps, non de quelque étudiant de première année inscrit à un cours d'introduction à la théorie des nombres. Même s'il n'est sûrement pas immunisé contre les erreurs, celle-ci semble plutôt extravagante. Peut-être que si nous creusons un peu plus profond, nous allons trouver qu'il voit les choses plus clairement que nous. Retournons en arrière et vérifions ce qu'il fait, à commencer par sa définition de $a:b=c:d$.

Dans notre analyse jusqu'à maintenant, nous avons considéré que $a:b=c:d$~signifie $a/b=c/d$, l'\emph{égalité des fractions} habituelle $ad=b$c. C'était tout à fait naïf. N'importe qui d'autre, y compris Zeuthen, se rend compte qu'Euclide avait deux notions très différentes de proportion, une dans le Livre~V qui traite de grandeurs arbitraires, et une dans le Livre~VII qui traite des nombres. Celle dans le Livre~V est la célèbre théorie grecque des proportions qui a été développée pour prendre en compte les grandeurs incommensurables. Cette théorie qui a des similitudes avec la théorie moderne des nombres réels, est usuellement associée au nom d'Eudoxe. Celle du Livre~VII traite des nombres, qui sont des grandeurs commensurables -- en fait, elles sont toutes multiples d'une grandeur unité fixée. Il a souvent été suggéré que la théorie du Livre~VII est plus ancienne, peut-être due aux pythagoriciens.

Euclide n'entend certainement pas $ad=bc$ quand il écrit «$a$~est à~$b$ comme $c$~est à~$d$» dans le Livre~VII. Il a défini ce que nous appellerons la \emph{proportionnalité pythagoricienne}, pour la distinguer de la \emph{proportionnalité eudoxienne} du Livre~V. Nous utilisons l'adjectif «pythagoricien» pour indiquer que cette proportionnalité a affaire uniquement à des nombres entiers. Voici la définition d'Euclide sous une forme moderne.
\begin{dvingt}[proportionnalité pythagoricienne]
  Nous disons que \(a:b=c:d\) s'il existe \(x\),~$y$, $m$ et~\(n\) tels que
  \begin{gather*}
%    \SwapAboveDisplaySkip
    a=mx,\quad b=nx,\\
    c=my,\quad d=ny.
  \end{gather*}
\end{dvingt}
\noindent Notez que nous pouvons comprendre cette définition comme disant que les fractions $a/b$~et~$c/d$ ont une simplification commune, à savoir~$m/n$. Clairement ceci implique que $ad=bc$ (égalité des fractions). La réciproque, bien qu'elle soit vraie pour les nombres entiers, est beaucoup plus profonde et échoue dans d'autres situations multiplicatives où la factorisation unique en nombres premiers n'est pas valable (voir les exemples 1~et~2). Observez aussi que la proportionnalité pythagoricienne dit simplement que $a$~est $m$~$n$ièmes parties de~$b$ et que $c$~est les mêmes parties de~$d$. C'est plus ou moins comment l'a en fait énoncé Euclide.

Euclide a besoin de l'équivalence de la proportionnalité pythagoricienne et de l'égalité des fractions dans la démonstration de son lemme. Rappelons comment il a déduit le lemme de la proposition~20. Il a supposé que $p$~était un nombre premier et que $ab=pc$, soit $a/p=c/b$ (égalité de fractions). Il a alors fait appel à la proposition~20. Mais la proposition~20 concerne la proportionnalité pythagoricienne, non l'égalité des fractions. La proposition dans les \emph{Éléments} qui nous assure que la proportionnalité pythagoricienne est la même que l'égalité des fractions est la suivante:
\begin{dixneuf}
  On a \(a:b=c:d\) si et seulement si \(ad=bc\).
\end{dixneuf}
\noindent Nous allons revenir sous peu à la proposition~19, mais d'abord voyons pourquoi la démonstration de la proposition~20 est correcte, en suivant Zeuthen. Pour le point crucial, nous savons que $u:c=v:d$, donc il existe $x$, $y$, $m$, et~$n$ tels que
\begin{gather*}
  u=mx,\quad c=nx,\\
  v=my,\quad d=ny.
\end{gather*}
Cela signifie que $n$~divise à la fois $c$~et~$d$ par définition! Nous n'avons pas besoin de le démontrer. En outre, $m>1$ parce que nous avons supposé, par contradiction, que $u$~ne divise pas~$c$. Naturellement, nous devons maintenant nous inquiéter de savoir si la proposition~13 (\emph{alternando}) est vraie, parce qu'elle ne dit pas simplement que si $a/b=c/d$, alors $a/c=b/d$. Cependant, on voit facilement que \emph{alternando} est vraie en échangeant les rôles de $m,n$~et~$x,y$ dans les équations de la définition~20.

Euclide n'a pas démontré \emph{alternando} de cette manière parce qu'il voyait les nombres $a$,~$b$, $c$, $d$, $x$ et~$y$ comme des grandeurs, quoique toutes multiples d'une même unité, tandis que $m$~et~$n$ étaient des choses qui répondaient à la question: combien? Ainsi le nombre~$x$ est une partie du nombre~$a$, et $m$~nous dit combien de~$x$ il faut pour faire~$a$. Des objets tels que $m$~et~$n$ ont été appelés dans la littérature moderne «nombres de répétition» par Fowler \cite{fowler99} et «scalaires» par Bashmakova \cite{bashmakova48}. Fowler appelle les nombres usuels des «nombres cardinaux». En fait, pour un étudiant d'aujourd'hui, il pourrait être utile pour comprendre cette distinction de penser aux nombres comme représentés par des \emph{ensembles} finis. Quand nous multiplions un ensemble~$a$ par un scalaire~$m$ nous prenons la réunion de $m$~copies disjointes de~$a$. Lorsqu'Euclide voulait multiplier deux nombres $a$~et~$b$, il considérait le nombre~$m$ d'unités présentes dans~$a$ et il posait $ab=mb$.

Comment alors Euclide démontre-t-il la proposition~19, que la proportionnalité pythagoricienne est équivalente à l'égalité des fractions? Puisque cette proposition est tout ce dont nous avons besoin pour achever la démonstration du lemme d'Euclide, certainement nous allons voir un appel au porisme ici, soit de manière directe, soit de manière indirecte. La moitié intéressante de la proposition est l'implication de $ad=bc$ à $a:b=c:d$. La démonstration d'Euclide procède comme suit. Si $ad=bc$, alors certainement $ac:ad=ac:bc$. Par ailleurs il est évident par la définition de la proportionnalité pythagoricienne que $ac:ad=c:d$ (prendre $x=ay$) et $ac:bc=a:b$. Par conséquent, $a:b=c:d$.

Stupéfiant! Aucun appel d'aucune sorte au porisme. Euclide nous a fourni une preuve de son lemme qui, ironiquement, ne dépend pas de manière essentielle de l'algorithme qui porte son nom, alors même qu'il commence le Livre~VII par cet algorithme et ses conséquences algébriques pour le plus grand commun diviseur. Pouvons-nous croire cela? Cela semble comme de la magie. Où le travail a-t-il été fait?

Bon, le symbolisme que nous avons adopté --~en utilisant le signe d'égalité pour une proportion~-- est trompeur et peut vous avoir compliqué la tâche de mettre le doigt sur l'erreur dans l'argument d'Euclide. Euclide lui-même disait, «Des choses égales à une même chose sont aussi égales l'une à une autre». Parce que nous avons internalisé cet axiome, il est dangereux d'utiliser un signe d'égalité dans une situation où la transitivité n'a pas lieu de manière évidente. Aussi, veuillez accepter nos excuses. En fait c'est une tâche non triviale de démontrer que la proportionnalité pythagoricienne est transitive. Essayez par vous-mêmes. C'est vrai pour les entiers naturels, mais à l'instar du lemme d'Euclide, c'est faux dans des situations multiplicatives plus générales, comme dans les deux exemples qui suivent. En outre, parce que l'égalité des fractions \emph{est} transitive, la proposition~19 est également en défaut.

Peut-être le contexte le plus simple dans lequel le théorème fondamental de l'arithmétique est en défaut est celui de notre premier exemple.

\begin{exemple}
  Considérons le monoïde multiplicatif des nombres entiers positifs 1,~4, 7, 10,~\dots\ qui sont congrus à~1 modulo~3. Dans ce monoïde, les nombres 4,~10 et~25 sont premiers, et $4\cdot25=10\cdot10$. La proportionnalité pythagoricienne n'est pas transitive: le lecteur peut vérifier que
  \[
  4:10=4\cdot25:10\cdot25=10\cdot10:10\cdot25=10:25,
  \]
  tandis que $4:10\neq10:25$ parce que 4,~10 et~25 sont premiers. En outre, les nombres 40~et~100 n'ont pas de pgcd algébrique. En fait, les diviseurs communs de 40~et~100 sont exactement 4~et~10, et aucun des deux ne divise l'autre.
\end{exemple}
Notre second exemple est plus compliqué, mais peut-être plus satisfaisant parce que c'est un système dans lequel vous pouvez aussi additionner.
\begin{exemple}
  Considérons le semi-anneau~$S$ des nombres réels~$a+b\sqrt2$ tels que $a$~et~$b$ sont des entiers positifs ou nuls, non tous deux nuls. Ici nous avons
  \[
  7(5+2\sqrt2)=(3+8\sqrt2)(1+2\sqrt2)
  \]
  et les quatre facteurs sont tous premiers. Notez que $S$~n'est pas l'anneau~$\mathbb Z[\sqrt2]$ d'entiers algébriques dans lequel $a$~et~$b$ peuvent être négatifs, anneau dans lequel le théorème fondamental de l'arithmétique est vérifié. Le seul élément inversible de~$S$ est~1, tandis que~$1+\sqrt2$ et toutes ses puissances sont inversibles dans~$\mathbb Z[\sqrt2]$. La proportionnalité pythagoricienne échoue dans~$S$ pour les mêmes raisons que dans l'exemple~1. En outre, $7(5+2\sqrt2)$~et~$7(1+2\sqrt2)$ n'ont pas de pgcd algébrique.
\end{exemple}

\section{Comment réparer l'argument d'Euclide.}

Il est intéressant qu'Euclide démontre explicitement, dans la proposition~11 du Livre~V, que la proportionnalité \emph{eudoxienne} est transitive mais qu'il omet de donner une démonstration que la proportionnalité \emph{pythagoricienne} est transitive, alors même que la démonstration de la proposition~19 fait un usage essentiel de ce fait. Y a-t-il un argument simple que nous pourrions incorporer au Livre~VII pour montrer que la proportionnalité pythagoricienne est transitive?

Bashmakova \cite{bashmakova48} pensait que c'était le cas. Elle a identifié l'utilisation non justifiée de la transitivité dans la démonstration de la proposition~19 comme étant le problème, contrairement à l'affirmation de Zeuthen selon laquelle la démonstration de la proposition~20 souffrait d'une pétition de principe. Elle a alors suggéré d'utiliser la transitivité de la proportionnalité eudoxienne pour rectifier l'erreur. À cette fin elle a donné une démonstration directe, qui n'invoque pas le porisme, que la proportionnalité pythagoricienne implique la proportionnalité eudoxienne. Mais si cette approche fonctionnait, nous aurions toujours une démonstration du lemme d'Euclide qui ne fait pas appel au porisme. Le problème est que Bashmakova n'a pas prouvé la réciproque, à savoir que la proportionnalité eudoxienne pour les nombres implique la proportionnalité pythagoricienne. C'est à cette partie de l'équivalence qu'elle devait réellement se confronter. En fait, il est facile de démontrer, sans utiliser le porisme, que la proportionnalité eudoxienne pour les nombres est équivalente à l'égalité des fractions, donc la proposition~19 peut être considérée comme disant que proportionnalités eudoxienne et pythagoricienne sont équivalentes. Dans cette perspective, l'erreur d'Euclide dans la démonstration de la proposition~19 se produit lorsqu'il démontre que la proportionnalité eudoxienne implique la proportionnalité pythagoricienne --~l'autre moitié de la démonstration est juste. En somme, la réparation de Bashmakova ne répare rien. (Nous avons découvert l'article de Bashmakova à partir d'une référence dans Narkiewicz \cite{narkiewicz00}).

Heath n'a fait aucun commentaire sur l'usage non fondé de la transitivité par Euclide. Cependant, dans ses notes sur la proportionnalité eudoxienne \cite[pages~126-129]{heath26}, il a donné une démonstration, attribuée à R. Simson, que la proportionnalité pythagoricienne est la même que la proportionnalité eudoxienne appliquée aux nombres entiers. Mais cette démonstration se trompe fatalement à la fin, lorsque les différentes définitions de \emph{partie} dans les Livres~V et~VII sont confondues: dans le Livre~V une \emph{partie} d'une grandeur est n'importe quel sous-multiple d'une autre grandeur, tandis que dans le Livre~VII une \emph{partie} d'un nombre doit être un autre nombre.

Les deux notions de proportionnalité ne se mettent pas facilement en relation, même si beaucoup d'auteurs ont été tentés d'imaginer que si. La proportionnalité pythagoricienne parle de divisibilité des nombres et est taillée pour étudier la factorisation. La proportionnalité eudoxienne, qui pour les nombres est équivalente à l'égalité des fractions, ne dit rien concernant la factorisation. Le fait qu'elles sont équivalentes pour les nombres est essentiellement le contenu de la proposition~19.

Quelle est la meilleure manière de réparer l'argument d'Euclide? Nous suggérons une manière qui utilise le porisme d'Euclide, que le plus grand commun diviseur est un plus grand commun diviseur algébrique, comme discuté dans la section~2. L'idée, que l'on peut trouver aussi dans \cite[Théorème~205]{taisbak71}, est de montrer que si un choix quelconque de $x$~et~$y$ établit la proportion pythagoricienne $a:b=c:d$, alors les choix canoniques~$x=\pgcd(a,b)$ et~$y=\pgcd(c,d)$ aussi. De cette manière, la transitivité a lieu, et la démonstration de la proposition~19 est réparée. En effet, la transitivité pouvait être en défaut uniquement si nous étions forcés d'utiliser un diviseur commun de $c$~et~$d$ dans la proportion~$a:b=c:d$ qui serait différent de celui utilisé dans la proportion~$c:d=e:f$. Si nous pouvons toujours utiliser le plus grand commun diviseur, alors la transitivité a clairement lieu.
\begin{fix}
  Supposons que \(a:b=c:d\). Si \(a=p\pgcd(a,b)\) et~\(b=q\pgcd(a,b)\), alors \(c=p\pgcd(c,d)\) et~\(d=q\pgcd(c,d)\).
\end{fix}
\begin{proof}
Par définition il existe $m$,~$n$, $x$ et~$y$ tels que
\begin{gather*}
  a=mx,\quad b=nx\text,\\
  c=my,\quad d=ny.
\end{gather*}
Parce que $x$~est un diviseur commun de $a$~et~$b$, et que $y$~est un diviseur commun de $c$~et~$d$, le porisme nous dit que nous pouvons trouver $i$~et~$j$ tels que $\pgcd(a,b)=ix$ et $\pgcd(c,d)=jy$. La première chose que nous voulons faire est de montrer que $i=j$. Par symétrie, il suffit de montrer que $i$~divise~$j$. Or $ix$~divise~$mx$, donc $iy$~divise~$my$. Pareillement, $ix$~divise~$nx$, donc $iy$~divise~$ny$. Ainsi $iy$~divise à la fois $my=c$~et~$ny=d$. Du porisme nous concluons que $iy$~divise~$\pgcd(c,d)=jy$, donc $i$~divise~$j$. Finalement, $c=p(iy)=p(jy)=p\pgcd(c,d)$ et $d=q(iy)=q(jy)=q\pgcd(c,d)$.
\end{proof}

Avec quel genre de démonstration du lemme d'Euclide nous retrouvons-nous? Voyons ce que nous avons fait. La proportionnalité pythagoricienne est essentiellement une relation entre fractions (pas entre nombres rationnels), \emph{a priori} plus forte que l'équivalence usuelle. Elle dit que $a/b$~est en relation avec~$c/d$ si $a/b$~et~$c/d$ ont une simplification commune, à savoir~$m/n$ avec $m$~et~$n$ comme dans la définition~20. (Naturellement, Euclide ne l'aurait jamais exprimé de cette manière, parce que pour lui $m$~et~$n$ sont des entités d'un type différent de $a$~et~$b$, comme cela a été discuté après la proposition~19.) Le porisme de l'algorithme d'Euclide peut être utilisé pour montrer que cette relation est transitive via notre réparation, et la transitivité est utilisée pour montrer qu'elle est équivalente à l'égalité des fractions $ad=bc$ (proposition~19). Alors, parce que la proportionnalité pythagoricienne est équivalente à l'égalité des fractions, nous voyons que si nous simplifions~$a/b$ autant que possible, nous obtenons la plus petite fraction équivalente à~$a/b$ pour l'égalité des fractions (ayant les plus petits termes possibles). Ceci n'est pas clair sans la proposition~19: il aurait pu se faire que $a$~et~$b$ soient premiers entre eux, et que cependant $a/b$~ne soient pas les plus petits termes possibles au sens qu'il y ait des nombres $c$~et~$d$ plus petits tels que $a/b=c/d$. (Dans l'exemple~1, les nombres 10~et~25 sont premiers entre eux, mais la fraction formelle~$10/25$ n'a pas les plus petits termes possibles parce que $10/25=4/10$.) L'équivalence des deux conditions (i)~que $a$~et~$b$ sont premiers entre eux et (ii)~que $a/b$~a les plus petits termes possibles est au cœur de la démonstration d'Euclide.

Pour établir le lemme d'Euclide, nous supposons que $p$~est un nombre premier qui divise~$ab$, disons $ab=pc$. Alors $a/p=c/b$, d'où il suit que $a/p$~et~$c/b$ ont une simplification commune, parce que l'égalité des fractions implique la proportionnalité pythagoricienne (proposition~19). Alors ou bien $p$~divise~$a$ et nous pouvons conclure, ou bien $p$~ne divise pas~$a$, et nous ne pouvons pas simplifier~$a/p$ parce que $p$~est premier. Dans ce dernier cas $c/b$~doit se simplifier en~$a/p$ (proposition~20), ce qui implique que $p$~divise~$b$.

\section{Parties canoniques.}

Insister sur le plus grand commun diviseur suggère que peut-être Euclide avait les parties canoniques toujours à l'esprit quand il utilisait la proportion pythagoricienne. (Ceci correspond plus ou moins à notre façon usuelle de penser à une fraction comme étant sous sa forme réduite.) Si l'on réclamait systématiquement des parties canoniques, alors la transitivité nécessaire dans la démonstration de la proposition~19 serait triviale, et il y aurait à peine besoin de la mentionner. En fait, Zeuthen (dans un article ultérieur \cite{zeuthen10}) et Itard \cite{itard61} ont proposé cette interprétation. Ils croyaient que quand Euclide montrait comment construire le plus grand commun diviseur en utilisant l'algorithme d'Euclide et qu'il donnait ses propriétés algébriques dans le porisme, il montrait simultanément comment interpréter la proportion $a:b=c:d$. En fait, dans la démonstration de la proposition~4 près du début du Livre~VII, Euclide a écrit~$a$ comme $m$~$n$ièmes parties de~$b$ en utilisant l'algorithme d'Euclide pour construire $b/n=\pgcd(a,b)$.

Pourquoi l'interprétation de la proportionnalité pythagoricienne en termes de parties canoniques ne résout-elle pas tout le problème? Bashmakova a envisagé cette idée mais l'a rejetée, en partie à cause d'une autre proposition d'Euclide dans le Livre~VII.
\begin{six}
  Si \(a:b=c:d\), alors \(a:b=(a+c):(b+d)\).
\end{six}
\noindent Cette proposition, qui se démontre facilement pour la proportionnalité pythagoricienne en suivant la démonstration d'Euclide, demande un développement additionnel substantiel pour une démonstration dans l'interprétation parties canoniques. Itard \cite{itard61} souligne ce problème. Bien que ceci ne ressorte pas de notre présentation, la proposition~6 est essentielle pour le développement d'Euclide dans la démonstration de son lemme. La démonstration que nous avons donnée de la proposition~13 (\emph{alternando}) n'est pas celle d'Euclide, et quoiqu'elle fonctionne pour la proportionnalité pythagoricienne, elle est inadéquate dans l'interprétation parties canoniques. La démonstration bien différente de \emph{alternando }par Euclide, qui est valable pour les deux interprétations (comme la plupart de ses propositions), se fonde sur la proposition~6.

En outre, \emph{alternando} est requis comme une étape clé dans la démonstration de la proposition~19 dans l'interprétation parties canoniques: l'étape où l'on constate que $ac:ad=c:d$. C'est clair pour la proportionnalité pythagoricienne mais pas pour les parties canoniques. Cependant, $ac:c=ad:d$ est clair pour les parties canoniques (les plus grands communs diviseurs sont $c$~et~$d$), et \emph{alternando} transforme cette affirmation en $ac:ad=c:d$. De fait, Euclide passe par \emph{alternando} pour démontrer la proposition~19. De cette manière, si sa démonstration de la proposition~6 est fautive, sa démonstration du lemme d'Euclide l'est en définitive de même.

Qu'est-ce qui est faux dans la démonstration de la proposition~6 dans l'interprétation parties canoniques? Supposons que
\begin{gather*}
%  \SwapAboveDisplaySkip
  a=mx,\quad b=nx,\\
  c=my,\quad d=ny,
\end{gather*}
\emph{où} $x=\pgcd(a,b)$ et $y=\pgcd(c,d)$. Alors clairement 
\begin{gather*}
  a=mx,\quad b=nx,\\
  a+c=m(x+y),\quad b+d=n(x+y),
\end{gather*}
mais nous devons encore vérifier que $x+y=\pgcd(a+c,b+d)$. Cela résulte de notre réparation, ainsi ce théorème sert également à réparer la démonstration d'Euclide de son lemme dans l'interprétation parties canoniques.

De cette manière le tableau d'ensemble des arguments d'Euclide qui mènent à la démonstration de son lemme est le suivant. La démonstration de la proposition~19 n'est pas valable pour la proportionnalité pythagoricienne, tandis que la démonstration de la proposition~6 n'est pas valable dans l'interprétation parties canoniques. Notre réparation, qui repose sur le porisme, établit que les deux interprétations sont équivalentes, sauvant de la sorte la ligne de raisonnement d'Euclide dans chaque interprétation. Ces deux interprétations, et notre réparation, ont été détaillées par Taisbak dans \cite{taisbak71}.

\section{Conclusion.}

En tentant de découvrir comment Euclide a démontré le lemme clé pour le théorème fondamental de l'arithmétique, nous avons été confrontés à la question de savoir si Euclide a utilisé, ou devait utiliser, l'algorithme d'Euclide d'une manière essentielle. Telle qu'elle est écrite, sa démonstration ne fait pas d'appel essentiel à l'algorithme. D'autre part, dans sa démonstration de la proposition~19, qui affirme que la proportionnalité pythagoricienne est équivalente à l'égalité des fractions (et donc à la proportionnalité eudoxienne), Euclide suppose sans le justifier que la proportionnalité pythagoricienne est transitive. Ce trou peut être comblé par un argument qui utilise le porisme de l'algorithme d'Euclide, et il semble raisonnable qu'Euclide aurait pu et aurait dû fournir un tel argument. Le fait que la proportionnalité pythagoricienne résulte de l'égalité des fractions a été appelé le \emph{Vierzahlensatz}\footnote{N.D.T. ``Théorème des quatre nombres''.}. Ce théorème est prouvé et développé en détail par Surányi \cite{suranyi90}, qui soutient qu'Euler l'avait noté.

La transitivité de la proportionnalité pythagoricienne a été mise en valeur par notre recherche. Son importance est mise en perspective par sa connexion avec deux autres propriétés multiplicatives: l'existence de pgcds algébriques et l'unicité de la décomposition en facteurs premiers. Dans un monoïde simplifiable commutatif vérifiant la condition de chaine des diviseurs\footnote{N.D.T. Il n'y a pas de suite indéfiniment décroissante pour la divisibilité.}, ces trois propriétés sont équivalentes. Pour les nombres naturels, d'autres approches complètement différentes du lemme d'Euclide et de l'unicité de la factorisation sont disponibles. On peut utiliser la récurrence, ou l'approche géométrique de Surányi \cite{suranyi90}.

Plusieurs commentateurs modernes n'ont pas vu le trou dans la transitivité dans la démonstration par Euclide de la proposition~19, ou le trou dans la proposition~6 dans l'interprétation alternative parties canoniques \cite{dijksterhuis30}, \cite{heath26}, \cite{knorr75}, \cite{mueller81}, \cite{vanderwaerden49}, \cite{zeuthen10}. Certains ont été conduits à croire qu'il était facile de démontrer que la proportionnalité pythagoricienne n'est qu'un cas particulier de la proportionnalité eudoxienne \cite{bashmakova48}, \cite{heath26}, \cite{knorr75} sans faire appel au porisme de l'algorithme d'Euclide. En fait, la proportionnalité pythagoricienne est \emph{a priori} plus contraignante que la proportionnalité eudoxienne appliquée aux nombres, comme nos deux exemples dans d'autres contextes multiplicatifs le montrent. Pour les entiers naturels, proportionnalités pythagoricienne et eudoxienne sont équivalentes, mais établir ce fait n'est pas chose triviale.

Comment Euclide peut-il avoir laissé un tel trou? Quand il a défini la proportionnalité eudoxienne pour les grandeurs, il a montré qu'elle était transitive (proposition~11 du Livre~V). Dans le Livre~VII, il a défini la proportionnalité pythagoricienne d'une manière complètement différente, mais il a supposé sans démonstration qu'elle aussi était transitive. Bien que la transitivité ne soit pas évidente, il aurait pu l'obtenir à partir du porisme avec lequel il a commencé le Livre~VII. Le fait qu'Euclide ait établi le porisme, mais qu'il ait échoué à l'utiliser au moment où il en avait le plus besoin, ne peut que laisser extraordinairement perplexe.

Finalement, il y a un consensus aujourd'hui que la proportionnalité eudoxienne est une idée sophistiquée qui a englobé la proportionnalité pythagoricienne plus simple et qu'elle l'a rendue obsolète. Notre analyse indique qu'au contraire, la proportionnalité pythagoricienne n'est pas un cas particulier immédiat de la proportionnalité eudoxienne. Elle est \emph{a priori} une relation strictement plus forte, particulièrement adaptée à l'étude de la divisibilité.

\begin{remerciements}
  Nous remercions Guram Bezhanishivili pour son aide dans la traduction du russe, Corné Kreemer pour son aide avec le néerlandais, et le personnel du prêt entre bibliothèques de nos universités pour leur soutien avec les ressources. Nous remercions également Andrzej Ehrenfeucht pour ses commentaires utiles.
\end{remerciements}

\bigskip

\noindent David Pengelley est professeur à l'Université d'État du Nouveau Mexique à Las Cruces. 

\noindent Fred Richman est professeur à l'Université Atlantique de Floride à Boca Raton.

\end{document}